# Back to harmonic mappings of compact Riemannian manifolds

Sergey E. Stepanov, Irina I. Tsyganok

**Abstract:** In this paper, we address several interconnected problems in the theory of harmonic maps between Riemannian manifolds. First, we present necessary background and establish one of the main results of the paper: a criterion characterizing when a smooth submersion or diffeomorphism between Riemannian manifolds is harmonic. This result provides a useful analytic condition for verifying the harmonicity of geometric mappings. Second, we investigate the $L^2$-orthogonal decomposition of the pullback metric associated with a harmonic map. We analyze the structure of this decomposition and discuss its geometric implications, particularly in the context of the energy density and trace conditions. Finally, we study harmonic symmetric bilinear forms and harmonic Riemannian metrics. Special attention is given to their role in the theory of harmonic identity maps. We derive new results that link these notions and demonstrate how they contribute to the broader understanding of harmonicity in geometric analysis.

## Introduction

Harmonic maps have emerged as a central concept in differential geometry, unifying and extending several classical ideas such as geodesics, minimal surfaces, and harmonic functions. Formally introduced in the seminal work of Eells and Sampson in 1964 [1], harmonic maps are defined as critical points of the energy functional associated with smooth mappings between Riemannian manifolds. Their natural variational origin endows them with deep geometric meaning and a wide range of applications (see [2] – [5]).

The theory of harmonic maps has since undergone significant development, becoming a rich area of research at the intersection of differential geometry, global analysis, and mathematical physics. Beyond their intrinsic geometric interest, harmonic maps play a fundamental role in the study of holomorphic maps between complex manifolds, stochastic analysis on manifolds, and non-linear sigma models in theoretical physics. These diverse connections underscore the broad mathematical relevance of the theory and motivate continued efforts to deepen our understanding of its structure.

In this paper, we revisit the geometry of harmonic maps from a new perspective, focusing on structural decompositions of geometric data associated with such maps. Our approach is inspired by classical ideas in global Riemannian geometry and aims to identify and clarify conditions under which certain classes of maps—particularly submersions, diffeomorphisms, and identity maps—are harmonic.

The structure of the paper is as follows. In Section 1, we provide background material and prove a general criterion for determining when a submersion or diffeomorphism between Riemannian manifolds is harmonic. Several examples of harmonic identity maps are also discussed to illustrate the criterion. Section 2 explores the geometric implications of the $L^2$-orthogonal decomposition of the pullback metric (the first fundamental form) induced by a harmonic map. Finally, in Sections 3 and 4, we introduce and analyze the notions of harmonic symmetric tensors and harmonic Riemannian metrics, with particular attention to their applications in the theory of harmonic identity maps.

The techniques developed here have the potential to contribute further to the understanding of harmonicity in geometric analysis and open new avenues for exploring geometric flows, stability questions, and variational problems in both mathematics and physics.

## 1. Harmonic maps

A systematic study of harmonic maps was initiated in 1964 by Eells and Sampson [1]. Detailed presentations of the results can be found in [2]; [3]; [4] and many other publications including monographs (see, for example, [5]). For definitions, notations, and results, we will refer to these works.

Let $(M, g)$ be a compact Riemannian manifold of dimension $n \geq 2$ with the Levi-Civita connection $\nabla$. Along with $(M, g)$, consider an $m$-dimensional ($m \geq 2$) Riemannian manifold $(\bar{M}, \bar{g})$ with the Levi-Civita connection $\bar{\nabla}$ and a $C^\infty$-mapping $f : (M, g) \to (\bar{M}, \bar{g})$. Let $f_*$ denote the differential of $f$ that is a $C^\infty$-section of the tensor bundle $T^*M \otimes f^*T\bar{M}$, where $f^*$ is the transpose map of $f_*$, and $f^*T\bar{M}$ the bundle with base $M$, fiber $T_{f(x)}\bar{M}$ over $x \in M$, and with the metric $\bar{g}'$ and connection $\bar{\nabla}'$ induced from $(\bar{M}, \bar{g})$.

The Riemannian structures $g$ and $\bar{g}$ define a metric $G$ on the fibers of the bundle $T^*M \otimes f^*T\bar{M}$. The connections $\nabla$ and $\bar{\nabla}'$ induce a connection $D$ in the bundle $T^*M \otimes f^*T\bar{M}$ such that (see [6] and [8])

$$(Df_*)(X,Y) = \bar{\nabla}'_X f_*Y - f_*\nabla_X Y$$

for any $X, Y \in C^\infty TM$. This formula defines the *second fundamental form* of the mapping $f$ (see [1] and [8]). In particular, a map with vanishing second fundamental form $Df_*$ is said to be *totally geodesic*. Such a map is characterized by the property that it carries geodesies to geodesies linearly (see [5, p. 9]). In turn, $g^* := f^*\bar{g}$ is called *the pullback metric* or, in other words, the *first fundamental form* of the mapping $f$ (see [1]). In particular, a map with vanishing first fundamental form is *constant* (see [1, pp. 112; 124]). The *energy density of the map $f$* is non-negative function given by $e(f) := trace_g g^*$ (see [8, p. 10]). In particular, if $e(f) = 0$, then $f$ is also a constant map.

The *energy of the mapping $f$* is determined by the integral formula

$$E(f) = \frac{1}{2}\int_M e(f)dv_g$$

**Remark 1.** The localization of this definition yields the Dirichlet energy functional $E_U(f)$, which is the energy of $f$ concentrated in a relatively compact open subset $U$ of $M$. The formula for $E_U(f)$ is obtained from (1.1) by replacing the integration over $M$ in (1.1) by the integration over $U$. For the functional $E_U(f)$ we have $0 \leq E_U(f) < \infty$ and $E_U(f) = 0$ if and only if $f$ is locally constant on $U$.

A smooth map $f : (M,g) \to (\bar{M}, \bar{g})$ of Riemannian manifolds is called *harmonic* if it gives an extremum of the energy functional $E_U(f)$ for any relatively compact open subset $U$ in $M$ with respect to all variations of $f$ with a compact support in $U$. The Euler-Lagrangian associated with the Dirichlet energy functional $E_U(f)$ is the *tension field $\tau(f) = trace_g(Df_*)$*. It is well-known that $f$ is called a harmonic map if and only if it satisfies the *Euler–Lagrange equation* $\operatorname{trace}_g(Df_*) = 0$. Let us prove the following statement.

**Theorem 1.** *Let $f: (M, g) \to (\bar{M}, \bar{g})$ be a smooth map from a compact Riemannian manifold $(M, g)$ of dimension $n$ and an m-dimensional complete Riemannian manifold $(\bar{M}, \bar{g})$ with dimension m. Suppose that f is harmonic map. Then the pullback metric satisfies the identity $\delta g^* = -\frac{1}{2} d\, e(f)$.*

*Conversely, if this identity holds and:*

- *$n = m$, map $f: (M, g) \to (\bar{M}, \bar{g})$ is a diffeomorphism, or*
- *$n \geq m$, map $f: (M, g) \to (\bar{M}, \bar{g})$ is a submersion,*

*then f is a harmonic map.*

**Proof.** The covariant derivative $\nabla g^*$ has the form (see [9, formulas (2.2) and (2.3)])

$$(\nabla_Z g^*)(X, Y) = \bar{g}'\big((Df_*)(Z, X), f_* Y\big) + \bar{g}'\big((Df_*)(Z, Y), f_* X\big)$$

for $g^* := f^* \bar{g}$ and any $X, Y, Z \in C^\infty TM$. It follows directly from this formula that

$$(\delta g^*)(Y) = -\bar{g}'(\tau(f), f_* Y) - \sum_{i=1}^n \bar{g}'\big((Df_*)(Y, X_i), f_* X_i\big);$$

$$\left(\nabla_Y (\operatorname{trace}_g g^*)\right) = 2 \sum_{i=1}^n \bar{g}'\big((Df_*)(Y, X_i), f_* X_i\big),$$

where $X_1, \ldots, X_n$ is an arbitrary local orthonormal basis of vector fields on $(M, g)$ and the operator $\delta: C^\infty(S^2 M) \to C^\infty(T^* M)$ is called the *divergence operator* on the vector space $S^2 M$ of symmetric bilinear differential two-forms. The operator $\delta$ is defined by the formula

$$(\delta \varphi)(Y) = -\sum_{i=1}^n (\nabla_{X_i} \varphi)(X_i, Y)$$

for an arbitrary $\varphi \in C^\infty(S^2 M)$. From the above equalities, we deduce the following equations:

$$(\delta g^*)(Y) + \frac{1}{2} \nabla_Y e(f) = -\bar{g}'(\tau(f), f_* Y). \tag{1.1}$$

since $e(f) = \operatorname{trace}_g g^*$. Then if $f$ is a harmonic, by virtue of (1.1), we have $\delta g^* = -\frac{1}{2} de(f)$ where $g^* = f^* \bar{g}$ is the pullback metric and $e(f) = \operatorname{trace}_g(f^* \bar{g}')$ is the energy density of the map $f$.

For a harmonic diffeomorphism $f: (M, g) \to (\bar{M}, \bar{g})$ for $n = m$ the converse of the first statement of our theorem is true (see [4]). Moreover, if $n \geq m$ and

$rang\,(f_*) \geq m$ at each point of $M$, then from the identities $\delta g^* + \frac{1}{2}de(f) = 0$ and (1.1) we obtain $\tau(f) = 0$. Therefore, the converse of the first statement of our theorem is also true if $n \geq m$ and $f : (M, g) \to (\bar{M}, \bar{g})$ is a harmonic submersion. The theorem is proved.

As an example of the above, consider the identity map $id_M : M \to M$, which is defined by the condition $id_M(x) = x$ for any point $x \in M$. In other words, it is a map that always returns to the same point that was used as its argument. It is obvious that $id_M : M \to M$ is a mapping that is both continuous and bijective. A differentiable bijection is not necessarily a diffeomorphism. However, the identity map $id_M : M \to M$ on a smooth manifold $M$ is indeed a diffeomorphism. Specifically, if we consider the identity map from a smooth manifold $M$ to itself - that is, to the same underlying set equipped with the same smooth structure - then this map is trivially smooth, bijective, and its inverse (which is itself) is also smooth. Therefore, it satisfies the definition of a diffeomorphism.

Furthermore, the identity map $id_M : (M, g) \to (M, g)$ is a trivial example of harmonic map. The cases of identity maps, there are close relations with infinitesimal harmonic transformations of $(M, g)$ (see Paragraph 2 of our paper).

The following problem arises: Given a compact Riemannian manifold $(M, g)$, does there exist a metric $\bar{g} \neq g$ on $M$ such that the identity map $id_M : (M, g) \to (M, \bar{g})$ is a harmonic map? The answers to this problem are given by the following two examples.

**Example 1.** The Ricci tensor has the following remarkable property, which follows directly from the reduced second Bianchi identity: $\delta Ric + \frac{1}{2}d\,s = 0$ for scalar curvature $s = trace_g\,Ric$ (see [9, p. 120]). Then, using the above theorem, we can conclude that the following statement is true: If $Ric$ is positive definite, then $\bar{g} = Ric$ is a metric $\bar{g}$ on $M$ and the identity map $id_M : (M, g) \to (M, \bar{g})$ is a harmonic map  Furthermore, if $Ric < 0$, then the result holds for the target manifold $(M, \bar{g})$, where $\bar{g} = -Ric$ (see also [10, pp. 86-87]).

**Example 2.** Now consider dimension 3. At any point, there exists an orthonormal frame $\{e_1, e_2, e_3\}$ such that the Riemann curvature tensor
$$Rm(e_i, e_j)e_k = \sigma(e_i, e_j)(\delta_{jk}e_i - \delta_{ik}e_j),$$
where $\delta_{ik} := g(e_i, e_k)$ and $\sigma(e_i, e_j)$ is the sectional curvature of the plane spanned by $e_i$ and $e_j$. In an (oriented) orthonormal frame $\{e_1, e_2, e_3\}$, the *cross curvature tensor* $Cr$ is defined by $Gr_{ij} = \left(\frac{detE}{\det g}\right)(E^{-1})_{ij}$ for the *Einstein tensor* $E = Ric - \frac{1}{2}s\,g$ (see [16, p. 88]). Then verify that $(Gr^{-1})^{ij}\nabla_{e_i}Gr_{jk} = \frac{1}{2}(Gr^{-1})^{ij}\nabla_{e_k}Gr_{ij}$. Using this, anyone can prove that if the sectional curvatures of $(M, g)$ are negative everywhere (or positive everywhere), then $Cr$ is a metric $\bar{g}$ on $M$ and $id_M : (M, \bar{g}) \to (M, g)$ is a harmonic map (see also [10, p. 88]).

## 2. The Ebin-Berger and York decompositions

Let $M$ be a smooth manifold. If $M$ is compact, the space $C^\infty(S^2M)$ of symmetric bilinear differential two-forms or, in other words, of two-covariant symmetric tensor fields on $M$, endowed with the $C^\infty$-topology, is a Fréchet space. Berger and Ebin proved in [11] a whole series of decompositions of that space and in particular they showed that, for any fixed Riemannian metric $g$ on a compact manifold $M$, the space $C^\infty(S^2M)$ splits into two orthogonal, complementary, closed subspaces (see also [9, p. 118]):
$$C^\infty(S^2M) = Im\,\delta^* \oplus Ker\,\delta \tag{2.1}$$
One of them is the image of the first-order differential operator $\delta^* : C^\infty(T^*M) \to C^\infty(S^2M)$ defined on the space $C^\infty(T^*M)$ of one-forms on $M$ by $\delta^*\theta := \frac{1}{2}L_\xi g$, where $L_\xi$ is the Lie derivative and $\xi = \theta^\#$ is the vector field dual (by $g$) to the 1-form (see [11]; [9, p. 117; 514]).

The other subspace is $Ker\,\delta := \delta^{-1}(0)$, where $\delta$ is the divergence operator $\delta : C^\infty(S^2M) \to C^\infty(T^*M)$; it is the adjoint of $\delta^*$ (see [11]; [12]) with respect to the usual $L^2$ inner product of symmetric differential two-forms defined by $g$ (see [9, p. 118])

$$\langle \varphi, \varphi' \rangle := \int_M g(\varphi, \varphi') dv_g,$$

where $\varphi, \varphi' \in C^\infty(S^2 M)$. Using the above we can prove the following corollary.

**Corollary 1.** *Let $f : (M, g) \to (\bar{M}, \bar{g})$ be a smooth map from a compact Riemannian manifold $(M, g)$ of dimension $n \geq 3$ and an m-dimensional complete Riemannian manifold $(\bar{M}, \bar{g})$ with dimension $m \geq 2$. Let $f^*\bar{g} = \delta^*\theta + \varphi_0$ be the Berger–Ebin $L^2$-orthogonal decomposition of the pullback metric $f^*\bar{g}$ on M. If $\xi = \theta^\#$ is an infinitesimal harmonic transformation, then the energy of such map f has the form $E(f) = C \, Vol(M, g)$ for a constant $C \geq 0$.*

**Proof.** Let $f : (M, g) \to (\bar{M}, \bar{g})$ be a harmonic map between a compact Riemannian manifold $(M, g)$ and a complete Riemannian manifold $(\bar{M}, \bar{g})$. From (1.3) we can deduce the $L^2$-orthogonal decomposition of the *pullback metric* $g^* = f^*\bar{g}$ as a symmetric bilinear form on $f^*T\bar{M}$:

$$g^* = \delta^*\theta + \varphi_0, \qquad (2.2)$$

where $\theta$ is some one-form $\theta \in C^\infty(T^*M)$ and $\varphi_0$ is some divergence free symmetric bilinear two-form $\varphi_0 \in C^\infty(S^2 M)$. Let $f : (M, g) \to (\bar{M}, \bar{g})$ be a harmonic map, then applying the operator $\delta$ to both sides of (2.2), we obtain

$$de(f) = -2\, \delta\delta^*\theta \qquad (2.3)$$

since $\delta g^* + \frac{1}{2} de(f) = 0$. At the same time, we recall that the *Sampson Laplacian* $\Delta_S : C^\infty(T^*M) \to C^\infty(T^*M)$ is defined by the identity (see [13])

$$\Delta_S := 2\delta\delta^* - d\delta.$$

Let us also recall that the vector field $\xi$ is an *infinitesimal harmonic transformation* in $(M, g)$ if the local one-parameter group infinitesimal transformations generated by $\xi$ is a group of harmonic transformations (see [13]). Furthermore, necessary and sufficient condition for a vector field $\xi$ on a Riemannian manifold to be infinitesimal harmonic transformation in $(M, g)$ is that $\theta \in Ker\, \Delta_S$, where $\theta$ is the g-dual one-form to $\xi$ (see [13]).

**Remark 1.** For example, the *Killing vector field* $\xi = \theta^\#$ is an infinitesimal harmonic transformation since $\theta \in Ker\, \Delta_S \cap Ker\, \delta$ (see [13]). We recall here that a vector field $\xi$ on $(M, g)$ is called a Killing vector field if the 1-parameter group of infinitesimal diffeomorphisms associated to $\xi$ consists of infinitesimal isometries (see [9, p. 40]). Moreover, all Killing vector fields on $(M, g)$ form a complete set of generators for the algebra which has a finite dimension. In conclusion, if $g^* \in Ker\ \delta$, then from (2.3) we obtain that $\xi = \theta^\#$ is the Killing vector field.

Using the above we can rewrite (2.3) as

$$de(f) = \Delta_S \theta + d\delta\theta.$$

Therefore, if $\xi$ is an infinitesimal harmonic transformation, then $e(f) = \delta\theta + C$ for a constant $C$. In this case $E(f) = C\ Vol(M, g)$ for $C \geq 0$ due to *Green's theorem* $\int_M \delta\theta\ dv_g = 0$ for a compact Riemannian manifold $(M, g)$. As a result, we can formulate the above corollary.

**Remark 2**. We recall that the extreme value of the extreme value of the energy of the harmonic immersion $f: (M, g) \to (\bar{M}, \bar{g})$ of compact manifold $(M, g)$ into a Riemannian manifold $(\bar{M}, \bar{g})$ has the form $E(f) = \frac{n}{2} Vol(M, g)$ (see [1, p. 124]).

In turn, *York's theorem* [15] is another well-known result and is also included in the monographs (see, e.g., [9, p. 130]). Namely, for any $n$-dimensional ($n \geq 3$) compact Riemannian manifold $(M, g)$ the decomposition holds (see [15] and [9, p. 130])

$$C^\infty(S^2 M) = (\text{Im }\delta^* + C^\infty M \cdot g) \oplus \left(\delta^{-1}(0) \cap \text{trace}_g^{-1}(0)\right) \quad (2.4)$$

where both terms on the right-hand side of (2.2) are infinite-dimensional and orthogonal to each other with respect to the $L^2$ inner scalar product. Furthermore, the second factor $\delta^{-1}(0) \cap trace_g^{-1}(0)$ of (2.5) is the vector space of *TT-tensors*.

**Remark 3**. We recall that a symmetric divergence free and traceless covariant two-tensor is called *TT*-tensor and is denote by $\varphi_{TT}$ (see, for example, [14]). The vector space of *TT*-tensors is defined by the condition

$$S_{TT}(M) := \{\varphi \in C^\infty(S^2 M) |\ \delta\ \varphi = 0, trace_g \varphi = 0\}.$$

As a consequence of a result of Bourguignon-Ebin-Marsden (see [9, p. 132]) the vector space $S_{TT}(M)$ is an infinite-dimensional vector space on any compact Riemannian manifold $(M, g)$. Such tensors are of fundamental importance in

General Relativity (see, for example, [14]; [15]; [16]) and in Riemannian geometry (see, for instance, [9, p. 346-347]).

Let $f : (M, g) \to (\bar{M}, \bar{g})$ be a harmonic map between a compact Riemannian manifold $(M, g)$ of dimension $n \geq 3$ and a complete Riemannian manifold $(\bar{M}, \bar{g})$. From (2.5) we can deduce the York $L^2$-orthogonal decomposition of the pullback metric $g^* = f^*\bar{g}$ as a symmetric bilinear form on $f^*T\bar{M}$:

$$g^* = (\delta^*\theta + \lambda\, g) + \varphi_{TT} \tag{2.5}$$

for some $\theta \in C^\infty(T^*M)$ and $\varphi_{TT} \in S_{TT}(M)$, and a scalar function $\lambda \in C^\infty(M)$. Applying the operator $trace_g$ to both sides of (2.5), we obtain

$$e(f) := trace_g g^* = -\delta\theta + n\,\lambda, \tag{2.6}$$

where $\theta$ is the $g$-dual one-form of $\xi$. Next, from (2.6) we obtain

$$\frac{1}{2} d(trace_g g^*) = -\frac{1}{2} d\delta\theta + \frac{n}{2}\, d\lambda. \tag{2.7}$$

In turn, applying the operator $\delta$ to both sides of (2.5), we obtain

$$\delta g^* = \delta\delta^*\theta - d\lambda. \tag{2.8}$$

Then from (2.7) and (2.8) we deduce the following equalities:

$$0 = \delta(f^*\bar{g}) + \frac{1}{2} d\left(trace_g f^*\bar{g}\right) = \frac{1}{2}(2\,\delta\delta^*\theta - d\delta\theta) - \frac{n-2}{2}\, d\lambda.$$

Therefore, we have the equation

$$\Delta_S \theta = (n-2)\, d\lambda. \tag{2.9}$$

If $n = 2$, then $\xi = \theta^\#$ is an infinitesimal harmonic transformation. In turn, if $n > 2$, then $\xi = \theta^\#$ is an infinitesimal harmonic transformation if and only if $\lambda = constant.$ In the case when $\lambda = constant,$ from (2.5) we obtain that $E(f) = \lambda\, Vol(M, g)$ for $\lambda \geq 0$ due to *Green's theorem* $\int_M \delta\theta\, dv_g = 0$ for a compact Riemannian manifold $(M, g)$. As a result, we formulate the corollary.

**Corollary 2.** *Let $f : (M, g) \to (\bar{M}, \bar{g})$ be a harmonic map between a compact Riemannian manifold $(M, g)$ of dimension $n \geq 2$ and an m-dimensional $(n \geq 2)$ complete Riemannian manifold $(\bar{M}, \bar{g})$. Let $g^* = (\delta^*\theta + \lambda\, g) + \varphi_{TT}$ be the York $L^2$-orthogonal decomposition of the pullback metric $g^* = f^*\bar{g}$ for some $\theta \in C^\infty(T^*M)$ and $\varphi_{TT} \in S_{TT}(M)$, and a scalar function $\lambda \in C^\infty(M)$. Then $\Delta_S \theta = (n-2)\, d\lambda$ for the Sampson Laplacian $\Delta_S$. Therefore, we conclude that*

- *if $n = 2$, then $\xi = \theta^{\#}$ is an infinitesimal harmonic transformation;*
- *if $n > 2$, then $\xi = \theta^{\#}$ is an infinitesimal harmonic transformation if and only if $\lambda = constant$;*
- *if $n > 2$ and $\xi = \theta^{\#}$ is an infinitesimal harmonic transformation, then the extreme value of the energy of the map $f$ has the form $E(f) = \lambda\, Vol(M, g)$ for the constant $\lambda \geq 0$.*

**Remark 4.** If $g^* = \varphi_{TT}$, then from (2.5) we obtain that $\delta^*\theta + \frac{1}{n}\delta\theta\, g = 0$. In this case $\xi = \theta^{\#}$ is a *conformal Killing vector field* (see [10, pp. 118-119]; [13]). We recall here that a vector field $\xi$ on $(M, g)$ is called a conformal Killing vector field if the one-parameter group of infinitesimal diffeomorphisms associated to $\xi$ consists of infinitesimal conformal transformations. On the other hand, if in formula (2.5) the second factor $\varphi_{TT} = 0$, then from (2.6) we obtain that $\Delta_S \theta = (n-2)\, d\lambda$.

In conclusion, we consider the identity harmonic map $id_M : (M, g) \to (M, \bar{g})$, where $\bar{g} = Ric$ for the positive definite Ricci tensor $Ric$. In this case $g^* = Ric$ and $e(id_M) = s$ for the scalar curvature $s = trace_g Ric$ of $(M, g)$. Then from (2.6) we obtain the equation $Ric = (\delta^*\theta + \lambda\, g) + \varphi_{TT}$, where $\Delta_S \theta = (n-2)\, d\lambda$. Therefore, if $n > 2$ and $\xi = \theta^{\#}$ is an infinitesimal harmonic transformation, then $\lambda$ is a constant and the *total scalar curvature* $s(M) := \int_M s\, dv_g$ (see [9, p. 119]) has the form $s(M) = \lambda\, Vol(M, g)$ for a constant $\lambda$ since $s = -\delta\theta + n\lambda$.

**Corollary 3**. *Let $id_M : (M, g) \to (M, Ric)$ be the identity harmonic map, where $(M, g)$ is a compact Riemannian manifold $(M, g)$ of dimension $n \geq 3$ and $Ric$ is the positive definite Ricci tensor $Ric$. Let $Ric = (\delta^*\theta + \lambda\, g) + \varphi_{TT}$ be the York $L^2$-orthogonal decomposition of the tensor Ricci $Ric$ for some $\theta \in C^\infty(T^*M)$ and $\varphi_{TT} \in S_{TT}(M)$, and a scalar function $\lambda \in C^\infty(M)$. If $\xi = \theta^{\#}$ is an infinitesimal harmonic transformation, then the energy $E(M)$ of the map $id_M : (M, g) \to (M, Ric)$ equals to the total scalar curvature $s(M)$ of $(M, g)$ and has the form $E(M) = n\lambda\, Vol(M, g)$, where $\lambda$ is constant.*

## 3. Harmonic tensors

A symmetric differential two-form $\varphi \in C^\infty(S^2M)$ is called *harmonic tensor* (see [12]), if $\delta\varphi = -\frac{1}{2}d\,(trace_g\,\varphi)$. It is easy to verify that the tensor field $\delta^*\theta$ is harmonic if and only if the $\xi = \theta^\#$ is an infinitesimal harmonic transformation of $(M,g)$. Namely, if $\varphi = 2\delta^*\theta$, then

$$\Delta_S\theta := 2\delta\delta^*\theta - d\delta\theta = 2\left(\delta\varphi + \frac{1}{2}d(trace_g\,\varphi)\right).$$

In this case the above statement is true. In addition, the Ricci tensor $Ric$ is another example of a harmonic tensor, since $\delta Ric + \frac{1}{2}d\,s = 0$ for the scalar curvature $s = trace_g\,Ric$ (see Example 1).

Let us define the vector space of harmonic tensors on $(M,g)$ by the condition (see [12])

$$\mathcal{H}_g = \left\{\varphi \in C^\infty(S^2M) \,\Big|\, \delta\varphi = -\frac{1}{2}d\,(trace_g\,\varphi)\right\}.$$

Along with the Ebin-Berger and decompositions of $C^\infty(S^2M)$, in this section we show another splitting of $C^\infty(S^2M)$. Space $\mathcal{H}_g$ will participate in this splitting. Let us consider the following operator: $\alpha_g: C^\infty(T^*M) \to C^\infty(S^2M)$, defined by $\alpha_g(\theta) = \delta^*\theta + \frac{1}{2}(\delta\theta)g$. It is easy to check that one can use the $L^2$ inner product to define the adjoint operator $\alpha_g^* = \delta + \frac{1}{2}d \circ trace_g$. Moreover, $Ker\,\alpha_g$ consists of infinitesimal isometrics (Killing vector fields). At the same time, $Ker\,\alpha_g^*$ consists of harmonic tensors. For example, $\alpha_g^*(Ric) = 0$.

**Remark 5.** The operator $B_g := -\alpha_g^*$ is a well-known operator of *Hamilton's Ricci flow theory* (see [7]; [10, p. 76-77]).

On the other hand, the principal symbol $\sigma_t(\alpha_{g_0})(e)$ of the operator $\alpha_{g_0}$ has the form $\sigma_t(\alpha_{g_0})(e) = \frac{1}{2}(t \otimes e + e \otimes t) + \frac{1}{2}(t \cdot e)g_0$ for all non-zero cotangent vectors $t$ and $e$ (see [11, p. 385]). It is easy to check that the principal symbol $\sigma_t(\alpha_{g_0})(e)$ is injective. It is well-known that (see [11, p. 383]) if $\alpha_g$ is a differential operator with injective symbol, then the $L^2$-orthogonal decomposition

$$C^{\infty}(S^2M) = Im\ \alpha_g \oplus Ker\ \alpha_g^* \qquad (3.1)$$

is true, where $\mathcal{H}_g = Ker\ \alpha_g^*$. There are some differences between operators with injective symbol and elliptic operators. An example is found in the following theorem (see [17]): Let $D$ be a differential operator with symbol that is injective but not surjective. Then $Ker\ D^*$ is infinite dimensional.

In our case $\mathcal{H}_g = Ker\ \alpha_g^*$ for the differential operator $\alpha_g$ with injective but not surjective symbol, therefore, $\mathcal{H}_g$ is infinite dimensional vector space.

**Remark 6.** In addition, we note that the space $S_{TT}(M) = \delta^{-1}(0) \cap trace_g^{-1}(0)$ is a subspace of $Ker\ \alpha_g^*$ and recall that it has infinite dimension (see [9, p. 132]) and hence $Ker\ \alpha_g^*$ must also have infinite dimensional.

As a result of the above, the following theorem will be true.

**Theorem 3.** *Let $(M, g)$ be a compact Riemannian manifold of dimension $n \geq 2$ and $\mathcal{H}_g$ be the vector space of harmonic tensor on $(M, g)$. Then the $L^2$-orthogonal decomposition $C^{\infty}(S^2M) = Im\ \alpha_g \oplus \mathcal{H}_g$ holds for the operator: $\alpha_g: C^{\infty}(T^*M) \to C^{\infty}(S^2M)$, defined by the identity $\alpha_g(\theta) = \delta^*\theta + \frac{1}{2}(\delta\theta)g$, where $\theta$ is an arbitrary one-form $\theta \in C^{\infty}(T^*M)$. Moreover, the vector space $\mathcal{H}_g$ is infinite-dimensional.*

From (3.1) we can deduce the $L^2$-orthogonal decomposition of the symmetric differential two-form $\varphi \in C^{\infty}(S^2M)$ as the following

$$\varphi = \left(\delta^*\theta + \frac{1}{2}(\delta\theta)g\right) + \varphi_h, \qquad (3.2)$$

where $\theta$ is some one-form $\theta \in C^{\infty}(T^*M)$ and $\varphi_h$ is some symmetric bilinear two-form $\varphi_h \in C^{\infty}(S^2M)$ such that $\delta\varphi_h = -\frac{1}{2}d(trace_g\ \varphi_h)$. In this case we have $\delta\varphi = \delta\delta^*\theta - \frac{1}{2}d(\delta\theta) + \delta\varphi^h = \frac{1}{2}\Delta_S\theta + \delta\varphi_h$ and hence $\frac{1}{2}\Delta_S\theta = \delta(\varphi - \varphi_h)$. Then $\xi = \theta^{\#}$ is an infinitesimal harmonic transformation if and only if $\delta\varphi = \delta\varphi_h$. Let $\varphi \in C^{\infty}(S^2M)$ be a *conservative tensor* which means $\delta\varphi = 0$ (compare with [12, p. 300]). Then $\Delta_S\theta = d(trace_g\ \varphi_h)$. As a result, $\xi = \theta^{\#}$ is an infinitesimal harmonic transformation if and only if $trace_g\ \varphi^h = C$ for a constant $C$. In this case from (3.2) we obtain $trace_g \varphi = -\frac{n-2}{2} \delta\theta + C$ and hence $\int_M (trace_g \varphi) dv_g = C\, Vol(M,g)$.

Let us consider the subbundle $\mathfrak{K}(M) \subset T^*M \otimes S^2M$ on a Riemannian manifold $(M,g)$, such that $K(X,Y,Z) = K(X,Z,Y)$ and $K_{12}(Z) = \sum_{k=1}^n K(X_k, X_k, Z) = 0$ for any $K \in \mathfrak{K}(M)$, arbitrary vector fields $X$, $Y$, $Z$ and a local orthonormal basis $\{X_1, \ldots, X_n\}$ of vector fields on $M$.

Using Weyl theorems on invariants (see [18]), we proved in [4] and [20] that $\mathfrak{K}(M)$ has pointwise irreducible with respect to the action of orthogonal group $O(n)$ decomposition

$$\mathfrak{K}(M) = \mathfrak{K}_1(M) \oplus \mathfrak{K}_2(M) \oplus \mathfrak{K}_3(M),$$

where

$$\mathfrak{K}_1(M) = \{K \in \mathfrak{K}(M) \mid K(X,Y,Z) + K(Y,Z,X) + K(Z,X,Y) = 0\},$$
$$\mathfrak{K}_2(M) = \{K \in \mathfrak{K}(M) \mid K(X,Y,Z) - K(Y,X,Z) = 0\},$$
$$\mathfrak{K}_3(M) = \{K \in \mathfrak{K}(M) \mid K(X,Y,Z) =$$
$$= \frac{1}{(n^2+n-2)}[(n+1)K_{12}(X)g(Y,Z) - K_{23}(Y)g(X,Z) - K_{23}(Z)g(X,Y)]\},$$
$$K_{23}(Z) = \sum_{k=1}^n K(Z, X_k, X_k) = 0 \tag{3.3}$$

Furthermore, using standard arguments of invariant theory one can prove that this the pointwise unique irreducible orthogonal decomposition of $\mathfrak{K}(M)$.

If we denote by $\bar\varphi = \varphi - \frac{1}{2}(trace_g \varphi) g$ for an arbitrary $\varphi \in \mathcal{H}_g$, then $\bar\varphi \in C^\infty(S^2M)$ and $\nabla\bar\varphi \in \mathfrak{K}(M)$. In this case, the covariant derivative $\nabla\bar\varphi$ is a cross-section of relevant invariant subbundles $\mathfrak{K}_1(M)$, $\mathfrak{K}_2(M)$ and $\mathfrak{K}_3(M)$, their direct sums $\mathfrak{K}_1(M) \oplus \mathfrak{K}_2(M)$, $\mathfrak{K}_1(M) \oplus \mathfrak{K}_3(M)$, $\mathfrak{K}_2(M) \oplus \mathfrak{K}_3(M)$.

We will say that the tensor field $\varphi \in \mathcal{H}_g$ belongs to the class $\mathfrak{K}_\alpha$ or $\mathfrak{K}_\alpha \oplus \mathfrak{K}_\beta$ for $\alpha, \beta = 1,2,3$ and $\alpha < \beta$ if $\nabla\left(\varphi - \frac{1}{2}(trace_g \varphi) g\right)$ is a cross-section of the invariant subbundles $\mathfrak{K}_\alpha(M)$ or $\mathfrak{K}_\alpha(M) \oplus \mathfrak{K}_\beta(M)$ for $\alpha, \beta = 1,2,3$ and $\alpha < \beta$, respectively. As a result, we invariantly define six classes of harmonic tensors.

We will give here interpretations of only three classes of harmonic tensors. In turn, descriptions of other classes of harmonic tensors can be found in our article [4].

Firstly, let $\varphi \in \mathfrak{K}_1$, then the covariant derivative of $\bar{\varphi} = \varphi - \frac{1}{2}(trace_g \varphi)g$ satisfies the equations $\delta^* \bar{\varphi} = \delta \bar{\varphi} = 0$. We recall here that the symmetric derivative $\delta^*: C^\infty S^2 M \to C^\infty S^3 M$ is defined by the equation $(\delta^* \omega)(X, Y, Z) := (\nabla_X \omega)(Y, Z) + (\nabla_Y \omega)(Z, X) + (\nabla_Z \omega)(Y, Z)$ for arbitrary $\omega \in C^\infty(S^2 M)$ and $X, Y, Z \in TM$ (see [9, p. 35; 356]). Therefore, $\bar{\varphi}$ is a *Killing tensor* with constant trace (see, for example, [24]). The geometry of manifolds bearing symmetric Killing tensor fields is described in detail in the literature (see, for example, [21] and [23]). Let us consider the *Sampson Laplacian* $\Delta_S: C^\infty(S^2 M) \to C^\infty(S^2 M)$, defined by the formula (see [24])

$$\Delta_S := \delta^* \delta + \delta \, \delta^*.$$

Then $\bar{\varphi} \in Ker\, \Delta_S$ for an arbitrary $\varphi \in \mathfrak{K}_1$ such that $\bar{\varphi} = \varphi - \frac{1}{2}(trace_g \varphi)g$. Therefore, the vector space $\mathfrak{K}_1$ is finite-dimensional (see [24, p. 23]).

We recall the following local result: if $(M, g)$ is a Riemannian manifold of constant curvature, then there exist a local coordinate system $x^1, \ldots, x^n$ in which the components $\varphi_{ij}$ of an arbitrary symmetric Killing tensor $\varphi \in C^\infty(S^2 M)$ can be expressed in the form (see [4])

$$\bar{\varphi}_{ij} = e^{2f}\left(A_{ijkl} x^k x^l + B_{ijk} x^k + C_{ij}\right)$$

where $f = \frac{1}{2(n+1)} ln(det\, g)$ and $A_{ijkl}$, $B_{ijk}$ and $C_{ij}$ are constants which satisfy the following identities:

$$A_{ijkl} = A_{jikl}, \quad A_{ijkl} = A_{ijlk}, \quad A_{ijkl} + A_{ikjl} = 0,$$
$$B_{ijk} = B_{jik}, \quad B_{ijk} + B_{ikj} = 0, \quad C_{ij} = C_{ji},$$
$$g^{ij} A_{ijkl} = g^{ij} B_{ijk} = 0.$$

for contravariant components $g^{ij}$ of $g$. Therefore, we can conclude that

$$\varphi_{ij} = \bar{\varphi}_{ij} - \frac{1}{n-2}(trace_g \bar{\varphi})\, g_{ij} =$$
$$= e^{2f}\left(A_{ijkl} x^k x^l + B_{ijk} x^k + C_{ij} - \frac{1}{n-2}(trace_g C)\, g_{ij}\right)$$

for the harmonic tensor $\varphi$ belonging to the class $\mathfrak{K}_1$.

Secondary, let $\varphi \in \mathfrak{K}_2$, then the covariant derivative of $\bar{\varphi} = \varphi - \frac{1}{2}(trace_g\varphi)g$ satisfies the equations $\delta\bar{\varphi} = 0$ and

$$(\nabla_X\bar{\varphi})(Y,Z) - (\nabla_Y\bar{\varphi})(X,Z) = 0 \qquad (3.4)$$

where $X, Y, Z$ are arbitrary vector fields on $M$. In this case, $\bar{\varphi}$ is a *Codazzi tensor* with constant trace (see [9, p. 435]). The geometry of Riemannian manifolds bearing Codazzi tensors is described in detail in the monograph [9].

Let us consider the *Bourguignon Laplacian* $\Delta_B: C^\infty(S^2M) \to C^\infty(S^2M)$, defined by the formula (see [9, p. 355])

$$\Delta_B := d\,\delta + \delta\,d$$

where $(d\,\omega)(X,Y,Z) := (\nabla_X\omega)(Y,Z) - (\nabla_Y\omega)(X,Z)$ for arbitrary $\omega \in C^\infty(S^2M)$ and $X, Y, Z \in TM$. Therefore, if the manifold $(M,g)$ is compact, then $\varphi \in \mathfrak{K}_2$ if and only if $\bar{\varphi} \in Ker\,\Delta_B$. In this case the vector space $\mathfrak{K}_2$ is finite-dimensional (see [9, p. 464]).

The following local result is well-known (see [9, p. 436]): an arbitrary Codazzi tensor $\bar{\varphi} \in C^\infty(S^2M)$ defined on a Riemannian manifold $(M,g)$ of constant curvature $C$ has the form $\bar{\varphi} = \nabla(dF) + C \cdot Fg$ for an arbitrary function $F \in C^\infty(M)$. In our case $F \in C^\infty(M)$ is a solution to the *Poisson equation* $\Delta F + n\,C \cdot F = const$ for the rough Laplacian $\Delta := div \circ grad$. Therefore, if $(M,g)$ is a Riemannian manifold of constant curvature $C$, then a harmonic tensor $\varphi$ of the class $\mathfrak{K}_2$ has the form

$$\varphi = \left(\nabla df - \frac{1}{n-2}(\Delta f)g\right) - \frac{2}{n-2}Cfg$$

for $\Delta f + n\,Cf = const$.

The third class of harmonic tensors is related to the following system of differential equations

$$\nabla_k\bar{\varphi}_{ij} = \frac{1}{2(n^2+n-2)}\left(n\,\nabla_k(trace_g\bar{\varphi})g_{ij} + (n-2)\nabla_j(trace_g\bar{\varphi})g_{ik} + (n-2)\nabla_i(trace_g\bar{\varphi})g_{kj}\right) \qquad (3.5)$$

for an arbitrary $\varphi \in \mathfrak{K}_3$. We introduce the symmetric tensor field $\omega$ with components

$$\omega_{ij} = \bar{\varphi}_{ij} - \frac{n}{(n^2+n-2)}(trace_g\bar{\varphi})g_{ij},$$

which, as easy follows from (3.5), satisfy the differential equations
$$2\nabla_k \omega_{ij} = \nabla_i(trace_g \omega) g_{kj} + \nabla_j(trace_g \omega) g_{ik}. \tag{3.6}$$
By well-known Sinyukov's theorem (see [22, p. 313; 330]), a Riemannian manifold $(M, g)$ bears a symmetric tensor field $\omega$ with components satisfying (3.6) if and only if $(M, g)$ admits a *projective diffeomorphism* to some (pseudo-)Riemannian manifold $(\widetilde{M}, \tilde{g})$. Furthermore, if a Riemannian manifold $(M, g)$ admits a projective diffeomorphism to some (pseudo-)Riemannian manifold $(\widetilde{M}, \tilde{g})$, then there exists a symmetric tensor field $\omega$ on $(M, g)$ the components of which satisfy equations (3.6). At the same time, we recall that a diffeomorphism $f: (M, g) \to (\widetilde{M}, \tilde{g})$ is called a projective mapping of $(M, g)$ onto $(\widetilde{M}, \tilde{g})$ if $f$ maps any geodesic curve in $(M, g)$ onto a geodesic curve in $(\widetilde{M}, \tilde{g})$.

## 3. Harmonic metrics

In the present section we are concerned with the space of harmonic Riemannian metrics on a compact $C^\infty$-manifold. If $\mathcal{M} \subseteq C^\infty(S^2 M)$ is the set of smooth Riemannian metrics on $M$ (those sections which at each point $x$ of $M$ induce a positive definite bilinear two-form on $T_x M$), it is well known that $\mathcal{M}$ is an open convex cone in $C^\infty(S^2 M)$.

In turn, Ebin's result from [19] asserts that if a smooth manifold $M$ is compact without boundary, then, at each metric $g_0 \in \mathcal{M}$, there is a slice for the usual action of the group $\mathcal{G}$ of diffeomorphisms of $M$ on the space $\mathcal{M}$. Then since $\mathcal{M}$ is an open convex cone in $C^\infty(S^2 M)$ it carries a natural structure of Fréchet manifold such that, for each fixed Riemannian metric $g_0 \in \mathcal{M}$ the tangent space at $g_0$, $T_{g_0}\mathcal{M}$, is the tangent space to the orbit $\mathcal{O}_{g_0}$ through $g_0$ is $T_{g_0}\mathcal{O}_{g_0} = Im\ \delta^*_{g_0}$ (see [19, p. 25]). So, the decomposition result can be read as $T_{g_0}\mathcal{M} = T_{g_0}\mathcal{O}_{g_0} \oplus (T_{g_0}\mathcal{O}_{g_0})^\perp$, where orthogonal is taken with respect to the $\mathcal{G}$-invariant metric $G$ given by
$$G_{g_0}(\varphi, \varphi') = \int_M trace(g_0^{-1}\varphi, g_0^{-1}\varphi')\, dv_{g_0},$$

Ebin's slices are of the form $exp_{g_0}(U)$ where $U$ is an open neighbourhood of zero in $(T_{g_0}\mathcal{O}_{g_0})^\perp$ and $exp_g$ is the exponential map, at $g_0$, of the metric $G$.

Let $g_0$ be a fixed Riemannian metric with Levi-Civita connection $\nabla^{g_0}$ on a compact manifold $M$. Another Riemannian metric $g$ on $(M, g_0)$ is called *harmonic with respect to* $g_0$ if the identity map $id_M : (M, g_0) \to (M, g)$ is a harmonic map (see [12, p. 296]). In accordance with the above, we denote by $\mathcal{M}_{g_0} \subset \mathcal{M}$ the vector space of harmonic metrics on $M$ with respect to $g_0$. Namely, we establish the following

$$\mathcal{M}_{g_0} = \left\{ g \in \mathcal{M} \,\middle|\, \delta_{g_0} g = -\frac{1}{2} d\left(trace_{g_0} g\right) \right\},$$

where $\delta_{g_0} := -trace_{g_0} \circ \nabla^{g_0}$ is the divergence operator of $(M, g_0)$. In this case $\mathcal{M}_{g_0} \subset \mathcal{H}_{g_0}$, where $\mathcal{H}_{g_0}$ is the infinite-dimensional vector space of harmonic tensors on $(M, g_0)$.

These questions remain open for investigation: What is the dimension of the moduli space $\mathcal{M}_{g_0}$ associated with a compact Riemannian manifold $(M, g_0)$? Additionally, how can one characterize (and compute) the dimension of the space of harmonic identity maps $id_M : (M, g_0) \to (M, g)$?

The answer to these problems is given in the example below.

Let $(M, g)$ be a Riemannian manifold and let $\varphi \in \mathcal{H}_{g_0}$, then $g(t) = g + t\varphi$ is a harmonic metric for any $t \in (-\varepsilon, \varepsilon)$ and for some $\varepsilon > 0$ provided that $(M, g)$ is compact (see [12, p. 297]). Then on a compact manifold $(M, g)$ a one-parameter smooth family of harmonic metrics $g(t)$ is defined for $t \in (-\varepsilon, \varepsilon)$. Therefore, $\mathcal{M}_g$ is the infinite-dimensional vector space on a compact manifold $(M, g)$. Furthermore, in this case $id_M(t) : (M, g) \to (M, g(t))$ for any $t \in (-\varepsilon, \varepsilon)$ is a one-parameter family of harmonic maps and $\mathcal{H}_g$ is infinite-dimensional. Then the following theorem is true.

**Theorem 4.** *Let $(M, g)$ be a compact Riemannian manifold and let $\mathcal{H}_g$ denote the infinite-dimensional vector space of tensors fields on $(M, g)$ that are harmonic with respect to the metric $g$. Then, for any fixed background metric $g$ and any $\varphi \in \mathcal{H}_g$,*

*there exists a one-parameter family of Riemannian metrics $g(t) = g + t\varphi$, defined for $t \in (-\varepsilon, \varepsilon)$, such that the corresponding one-parameter family of identity maps $id_M(t): (M, g) \to (M, g(t))$ forms an infinite-dimensional space of harmonic maps.*

From the above, it follows that a compact Riemannian manifold $(M, g)$ admits a one-parameter family of harmonic Riemannian metrics of the form $g(t) = g + t\,\varphi$, defined for all $t \in (-\varepsilon, \varepsilon)$, where $\varphi$ is a symmetric bilinear form belonging to the infinite-dimensional vector space $\mathcal{H}_g$ of tensors that are harmonic with respect to the metric $g$. Motivated by this observation, we consider variations of the background metric $g$. Unless otherwise specified, all covariant derivatives, inner products, and traces are taken with respect to $g$.

The first variation of the Ricci tensor $Ric$ has the form (see [9, p. 64])

$$-2\dot{Ric}_\varphi = \Delta_L \varphi_{ij} + \nabla_i \nabla_j (g^{kl}\varphi_{kl}) + \nabla_j(\delta\varphi)_i + \nabla_i(\delta\varphi)_j,$$

where $\Delta_L$ is the *Lichnerowicz Laplacian* defined on symmetric 2-tensors (see also [10, p. 69]), that is, the Sampson Laplacian (see [24])

$$\Delta_S \varphi_{ij} := \Delta \varphi_{ij} + 2g^{kl} R^m_{kij} \varphi_{ml} - g^{kl} R_{il}\varphi_{jk} - g^{kl} R_{jl}\varphi_{ik}.$$

Let $\varphi$ with local components $\varphi_{jl}$ be a harmonic tensor, then the first variation of the Ricci tensor $Ric$ can be rewritten in the form

$$-2\dot{Ric}_\varphi = \Delta_S \varphi_{ij}.$$

The first variation of scalar curvature has the form (see [9, p. 63])

$$\dot{s}_\varphi = -\Delta\left(trace_g \varphi\right) + \delta(\delta\varphi) - g(\varphi, Ric),$$

where $g(\varphi, Ric) = g^{ik} g^{jl} \varphi_{ij} R_{kl}$. Let $\varphi(x)$ be a harmonic tensor, then the first variation of scalar curvature can be rewritten in the form

$$\dot{s}_\varphi = -\frac{1}{2}\Delta\left(trace_g \varphi\right) - g(\varphi, Ric).$$

Furthermore, $(\dot{dv})_\varphi = \frac{1}{2}\left(trace_g \varphi\right) dv_g$ (see [9, p. 65]).

Next, we consider the first variation of the Einstein-Hilbert functional $\mathbf{S}_g = \int_M s\, dv_g$ under the variation of metric: $g(t) = g + t\,\varphi$, defined for all $t \in (-\varepsilon, \varepsilon)$, where $\varphi$ belongs to the infinite-dimensional vector space of harmonic tensors on $(M, g)$. We know that (see [9, p. 120]):

$$\dot{\mathbf{S}}_\varphi = \int_M \left(\dot{s}_\varphi + (s/2)\, trace_g \varphi\right) dv_g,$$

where in our case $\dot{s}_\varphi = -\frac{1}{2} \Delta\left(trace_g \varphi\right) - g(\varphi, Ric)$. Then using the divergence theorem, we obtain the following formula:

$$\dot{\mathbf{S}}_\varphi = -\int_M g(E_g, \varphi) dv_g$$

for the *Einstein tensor* $E_g = Ric - \frac{s}{2} g$. Then, if $(M, g)$ is compact, we can conclude that the critical metrics for $\mathbf{S}_g$ (under variations $\varphi \in \mathcal{H}_g$) are those with vanishing Einstein tensor.

**Corresponding Author: Sergey E. Stepanov:**
Department of Mathematics, Russian Institute for
Scientific and Technical Information of the
Russian Academy of Sciences, 125190 Moscow,
Russia, E-mail address: s.e.stepanov@mail.ru

**Irina I. Tsyganok:** Department of Mathematics,
Finance University, 125468 Moscow, Russia,
E-mail address: i.i.tsyganok@mail.ru